\documentclass[12pt]{amsart}
\usepackage{amsmath,euscript,amssymb,amsthm,amscd,graphicx,epsfig}
\usepackage{amsmath,euscript,amssymb,amsthm,xypic,amscd,graphicx,epsfig}
\newtheorem{theorem}{Theorem}[section]
\newtheorem{proposition}[theorem]{Proposition}
\newtheorem{lemma}[theorem]{Lemma}
\newtheorem{corollary}[theorem]{Corollary}
\theoremstyle{definition}
\newtheorem{remark}[theorem]{Remark}
\newtheorem{definition}[theorem]{Definition}

\newtheorem{example}[theorem]{Example}

\numberwithin{equation}{section}
\newcounter{llistadepth}

\newenvironment{manlist}[1]{\addtocounter{llistadepth}{1}
      \edef\llistacontador{llista\romannumeral\the\value{llistadepth}}
      \list{({#1{\llistacontador}})}{\usecounter{\llistacontador}
      \def\makelabel##1{\hss\llap{##1}}
      \itemsep=2pt\parsep=0pt\topsep=3pt plus 1pt minus 1 pt}}{\endlist
      \addtocounter{llistadepth}{-1}}

\newcommand{\bth}{\begin{theorem}}
\renewcommand{\eth}{\end{theorem}}
\newcommand{\bpr}{\begin{proposition}}
\newcommand{\epr}{\end{proposition}}
\newcommand{\bco}{\begin{corollary}}
\newcommand{\eco}{\end{corollary}}
\newcommand{\ble}{\begin{lemma}}
\newcommand{\ele}{\end{lemma}}
\newcommand{\bre}{\begin{remark}\rm}
\newcommand{\ere}{\end{remark}}
\newcommand{\bex}{\begin{example}\rm}
\newcommand{\eex}{\end{example}}
\newcommand{\bde}{\begin{definition}}
\newcommand{\ede}{\end{definition}}

\newcommand{\Q}{\mathfrak{Q}}
\renewcommand{\ss}{(S^1,S^n)}
\newcommand{\rr}{(\R,\R^n)}
\def\la#1{\hbox to #1pc{\leftarrowfill}}
\def\ra#1{\hbox to #1pc{\rightarrowfill}}
\def\fract#1#2{\raise3pt\hbox{$ #1 \atop #2 $}}

\newcommand{\im}{\mathrm{Imm}}
\newcommand{\imm}{\mathrm{Imm'}}
\newcommand{\emb}{\mathrm{Emb}}

\newcommand{\so}{SO(n)}

\renewcommand{\H}{\mathbb{H}}

\newcommand{\R}{\mathbb{R}}

\numberwithin{theorem}{section}

\title{On the  Homology of the Space of Singular Knots}
\author{Hossein Abbaspour \and David Chataur  }
\thanks{The second author is partially supported by ANR grant 06-JCJC-0042 ``Op\'erades, Big\`ebres et
Th\'eories d'Homotopie''.}
\address{Laboratoire Jean Leray, Universit\'e de Nantes, Nantes 44300, France.}
\email{abbaspour@univ-nantes.fr}

\address{Laboratoire Paul Painlev\'e,
Universit\'e de Lille 1, 59655 Villeneuve d'Ascq C\'edex, France}

\email{David.Chataur@math.univ-lille1.fr}

\begin{document} \maketitle
\begin{abstract}In this paper we introduce various associative products on the homology of the space of knots and singular knots in $S^n$. We prove that these products are related through a desingularization map. We also compute some of these products and prove the nontriviality of the desingularization morphism.
\end{abstract}

\tableofcontents

\section{Introduction}
Various authors have tried to introduce  a general framework for the structures introduced by Chas-Sullivan \cite{cs} and Cohen-Jones \cite{CJ} on the homology and on the equivariant homology of the free loop space of a closed manifold. One approach, first used in \cite{ks} and then formulated in  \cite{gs}, is the notion of fiberwise (homotopy) monoid. We use this formulation to introduce and calculate various algebra structures on the homology of knots, immersions, and singular knots with $k$ double points.

In this paper the \emph{knot space} $\emb(S^1,S^n)$ is the space of all embeddings of $S^1$ in the $n$-sphere $S^n$, and similarly $\im(S^1,S^n)$ is the space of all \emph{immersions} of $S^1$ in $S^n$.  For a knot or an immersion $\gamma:S^1\rightarrow S^{n}$, the point $\gamma(1)$ is called the \emph{marked point} of $\gamma$. Here $S^n$ is the unit sphere in $R^{n+1}$ and we make the identification $\R^n\simeq \{0\}\times \R^n \subset 
\R^{n+1}$.  We think of $SO(n)\subset SO(n+1)$  as the stabilizer of $(1,0,\cdots,0)$ and $SO(n-1)\subset SO(n+1)$ as the stabilizer of $(1,0,\cdots,0)$ and $(0,1,0,\cdots0)$.
$US^n$ is the unit tangent sphere bundle of $S^n$ which can be identified with $SO(n+1)/SO(n-1)$. Let $LUS^n=C^\infty (S^1, US^n)$ be the free loop space of $S^n$.

The subspace  $\imm  \ss \subset \im \ss$ consists of those immersions which are not singular at the marked point and $\im_{k}\ss\subset\imm \ss$ is the space of immersions with exactly $k$ double \emph{points}. 

We will be considering \emph{uncompactified} or \emph{long} version of the objects defined above which roughly corresponds to the one dimensional object obtained in $ \R^n$ using the \emph{stereographic} map with the marked point as the center point.

For instance, the space of \emph{long knots}\footnote{There are some other models for the space of long knots, all homotopy equivalent to ours.} in $\R^n$ is 
$$
\emb^{l} \rr=\{f:\R\times \{0\} \rightarrow\R^n, supp(f)\subset (-1,1)
\text{ and } f([-1,1]\times \{0\}^{n-1}
)\subset B^n\}
$$
where $supp(f)=\{t\in \R \subset \R^n| f(t)\neq (t,0,\cdots,0) \}$ and $B^n$ is the unit ball. Similarly one defines $\im^{l}\rr$ the space of long immersions and $\im^{l}_{k}\rr$ the space of long immersion with $k$ singular double points. 

The purpose of this paper is to introduce various associative algebra structures on the homology of $\emb \ss$, $\imm \ss$ and $\im \ss $, and finally singular knot with double singular points $\coprod _{k}\imm_k \ss$. In the case of singular knots, one naturally expects a compatibility with the Vassiliev spectral sequences. However this issue is not addressed in this paper.  In order to introduce these products, for each space we consider a homotopically equivalent model as the total space of a fibration whose fibres form a continous family of homotopy associative monoids and the base space is a closed oriented manifold. This should be compared with the loop space fibration  $\Omega M \rightarrow LM \rightarrow M$ of a closed oriented manifold.

\setcounter{section}{3} \setcounter{theorem}{0}
\bth

$(H_{*+2n-1} (\emb(S^1,S^n)), \mu_{em}(-,-) )$ is a graded commutative  algebra.
 \eth

\bth
$(H_{*+2n-1}(\imm (S^1,S^n)), \mu_{im}(-,-) )$  is  a graded commutative algebra. 
\eth

The commutativity of the product follows from the computation of the homology by spectral sequences arguments.

\setcounter{section}{4} \setcounter{theorem}{1}
\bth
\label{thm3}
$H_{*+2n-1}(\im (S^1,S^n))$  is a graded commutative algebra and the map induced by inclusion
$i_{*}:H_{*+2n-1}(\imm (S^1,S^n))\rightarrow H_{*+2n-1}(\im (S^1,S^n))$ is a map of algebras.
\eth

After restricting to knots with double singular points,  one obtains:
\bth
There is a collection of maps $$\mu_{im}^{k,l}: H_*(\imm_k \ss)\times  H_*(\imm_l \ss) \rightarrow  H_{*-2n+1}(\imm_{k+l} \ss)$$

such that 

$$
\mu_{im}^{k+l,m}(\mu_{im}^{k,l}(a,b),c))=\mu_{im}^{k,l+m}(a,\mu_{im}^{l,m}(b,c)).
$$

Said slightly differently, $(\H_{*}(\emb^{s}(S^1,S^n)),\mu^{s})$ is an associative algebra where

$$
\H_{*}(\emb^{s}(S^1,S^n))=\oplus_{k} H_{*+2n-1}(\imm_k \ss) \text{ and }  \mu^{s}= \oplus_{k,l} \mu_{im}^{k,l}.
$$
\eth
It is not obvious  why this algebra is or should be commutative.

In order to compare the algebra structure on the homologies of knot and singular knot  spaces we introduce a desingularization map. Informally speaking, we resolve a singular knot at a double point in all possible ways, parameterized by unit vectors perpendicular to the tangent plane at the singularity. Of course this map has a certain degree and it is not a map of algebras. However it verifies a certain compatibility condition with respect to the number of singularities and the product. More precisely,

\setcounter{section}{5} \setcounter{theorem}{5}
\bth
The desingularization morphisms $\sigma_k: H_*(\imm_k \ss) \rightarrow H_{*-k(n-3)}(\emb \ss)$, $k\geq 0$, are compatible with the products i.e. for $x\in H_*(\imm_k \ss)$ and $y\in H_*(\imm_{l}(\ss))$
$$
\mu_{em}(\sigma_k (x), \sigma_{l}(y))= \sigma_{k+l}(\mu_{im}^{k,l}.(x,y)).
$$
In other words, 
$$\sigma=\oplus  \sigma_{k}:\H_{*}(\emb^{s}(S^1,S^n))\rightarrow \H_{*}(\emb \ss),$$
is a map of algebras, where $\H_{*} (\emb \ss)=H_{*+2n-1 }(\emb \ss)$.

\eth

\begin{table}[htdp]
\begin{center}\begin{tabular}{|c|c|}
\hline $ \emb^l \rr $ & Long knots in $\R^n$ \\
 \hline $\emb(S^1,S^n)$ &  Embeddings in $S^n$ (Knots) \\ 
 \hline $\im^l \rr$ & Long immersions in $\R^n$ \\ 
 \hline $\imm \ss$ & Immersions in $S^n$ without singularity at $t=0$ \\ 
 \hline  $\im^l_{k} \rr$ & Long immersions in $\R^n$ with $k$ double points \\
 \hline $\imm_{k}\ss$ & Immersions in $S^n$ with $k$ double points away from $ t=0$ \\ 
 \hline 
 
  \end{tabular} 
  \bigskip
 
  \caption{Notation}
\end{center}

\label{Notation}
\end{table}

\textbf{Acknowledgments:} 
The first author would like to thank the Hausdorff Research Institute and the MPI in Bonn for their support while this paper was written up. The second author warmly thanks the Laboratoire Jean Leray at the University of Nantes for their invitation throught the Matpyl program. 

\setcounter{section}{1} \setcounter{theorem}{0}
\section{Long immersions and long knots}

We start this section with a few remarks on the homotopy type of the space of knots $\emb \ss$ and immersions $\im' \ss$, in terms of the long knots $\emb ^l \rr$ and long immersions $\im^l \rr$.
The main references are \cite{budney1, budney2, BC}. 

 Note that $SO(n-1)$ acts on $SO(n+1)\times \emb^l\rr$, by the natural action on the first factor.  It acts on the second factor via the identification $\R \times \R^{n-1}=\R^n \simeq  \{0\}\times \R^n \subset \R^{n+1}$.  One may form the quotient $SO(n+1)\times_{SO(n-1)} \emb^l \rr$ which has the same homotopy type of $\emb \ss$ (see \cite{budney2}).  A homotpy equivalence  is given by the stereographic map with
 $(1,0,\cdots,0)\in S^n \subset \R^{n+1}=\R\times \R^{n}$ as the center.  Given a pair $(A,\gamma) \in SO(n+1)\times_{SO(n-1)} \emb^l \rr$ one can compactify $\gamma$ in $\R^{n+1}$ to get a knot in $S^n\subset \R^{n+1}$ passing through $(1,0,\cdots, 0)$. Then, we let $A$ act on  $\gamma$ to get an element of $\emb\ss$. This defines a homotopy equivalence map 
 \begin{equation}
 \rho_{emb}:  SO(n+1)\times_{SO(n-1)} \emb^l \rr \rightarrow  \emb \ss
\end{equation}
which makes the following diagram commute.
 \begin{equation}
\xymatrix{SO(n+1)\times_{SO(n-1)} \emb^l \rr \ar[r]^-{ \rho_{emb}} \ar[d] & \emb \ss \ar[d]^-{ev} \\  SO(n+1)/ SO(n-1) \ar[r]  &US^n}
\end{equation}
Here $ev(\gamma)= (\gamma(0), \gamma'(0)/\| \gamma'(0)\|)$ for $\gamma \in \emb \ss$. Therefore, we have the isomorphism,
 \begin{equation}
 \rho_{em}: H_{*}(SO(n+1)\times_{SO(n-1)} \emb^l \rr ) \rightarrow   H_{*}(\emb\ss).
 \end{equation}
 
  Similarly, one can consider the Borel construction $SO(n+1)\times_{SO(n-1)} \im^l \rr$ and $SO(n+1)\times_{SO(n-1)} \im^l_{k} \rr$ the homotopy equivalences induced by the stereographic map which gives rise to the commutative diagram
  
\begin{equation}
\xymatrix{SO(n+1)\times_{SO(n-1)} \imm^l \rr \ar[r]^-{ \rho_{im}} \ar[d] & \imm \ss \ar[d]^-{ev} \\  SO(n+1)/ SO(n-1) \ar[r]  &US^n}
\end{equation}
  
Similarly, one has the isomorphisms,
  \begin{equation}
 \rho_{im}:  H_{*}(SO(n+1)\times_{SO(n-1)} \im^l \rr )  \rightarrow H_{*}(\imm\ss)
 \end{equation}
  \begin{equation}
   \rho^{k}_{im}: H_{*}(SO(n+1)\times_{SO(n-1)} \im^l_{k} \rr)   \rightarrow H_{*}(\imm_{k}\ss)\end{equation}\
 
\section{Stringy operations}

In this section we recall the general framework where some of the string topology operations can be defined. This was observed by several authors independently, including the authors of this article. However the main published reference is \cite{gs}, therefore we follow their approach.
Let $F\rightarrow E\overset{\pi}{\rightarrow} M$ be a fiber bundle over a closed compact manifold $M$ of dimension $d$. We suppose that $E$ is equipped with a fiberwise associative product, that is a fiber bundle map  $m:E\times_B E \rightarrow E$ such that 

\begin{equation}
m(x,m(y,z))=m(m(x,y),z)
\end{equation}

Since the associativity of the product is not achieved for most of the interesting cases, one must adapt this definition to the case where the product is homotopy associative and consider \emph{fiberwise homotpy monoids}. 

One can take one step further by acquiring an operadic approach in the monoidal category of spaces fibered over $M$. This is simply done by introducing the trivial bundles $\bar{\mathcal{C}}_n=\mathcal{C}_n \times M\rightarrow M$ where $\mathcal{C}_n$ is the space of $n$ little cubes (see \cite{gs}). Then, one  may consider a fiber bundle  $E\rightarrow M$ which is an algebra over the bundle operad $\bar{\mathcal{C}}=\coprod_n \bar{\mathcal{C}}_n$  \emph{i.e.} we have a collection of fibre bundle maps $ E\times \bar{\mathcal{C}}_n \rightarrow  E$ subject to the usual axioms for the algebras over an operad.

This formulation allows us to define a commutative and associative product of degree $-d$

$$
\bullet : H_*(E)\times H_*(E)\rightarrow H_{*-d}(E)
$$
and a bracket of degree $(1-d)$

$$
\{-,-\}:H_*(E)\times H_*(E)\rightarrow H_{*+1-d}(E)
$$
which satisfies the Jacobi identity for a Lie bracket of degree $1-d$. One could be more ambitious and consider algebras over the bundle operad of framed little cube operad in order to equip $H_*(E)$ with a BV-algebra structure. This feature is not discussed in this paper.

Applying this machinery to the fibration $$\emb  \ss \simeq SO(n+1)\times_{SO(n-1)} \emb^l \rr\rightarrow  SO(n+1)/SO(n-1),$$ 
one obtains immediately  a commutative and associative product of degree $-2n+1$
\begin{equation}
\begin{split}
\mu_{em}(-,-)= \rho_{em}( \rho_{em}^{-1}(-)\bullet \rho_{em}^{-1}(-)):&
 H_*(\emb\ss) \times H_*(\emb\ss)\rightarrow \\ &H_{*-2n+1}(\emb  \ss)
 \end{split}
\end{equation}

\begin{theorem}(Gruher-Salvatore, ...)
$( H_*(\emb\ss), \mu_{em})$  is a  commutative and associative algebra.
\end{theorem}

This product is also introduced  in Gruher and Salvatore's paper for the case $n=3$ based on the R. Budney's work \cite{budney1}. Note that for all $n\geq 3$, the little cube operad acts on the space of long knots $\emb^l \rr$, therefore this product is well-defined for all $n\geq 3$. 

Similarly, one can consider the fibration
$$
\imm\ss \simeq SO(n+1)\times_{SO(n-1)} \im^l \rr \rightarrow SO(n+1)/SO(n-1).
$$
which accommodates  a fiberwise homotopy associative product but not the action of the full little cube bundle operad.  For instance, we lose the homotopy commutativity of the product. In the case of embeddings, the product at  the level of fibers is basically the connected sum of long knots. 
 
\begin{equation}
\begin{split}\mu_{im}(-,-)= \rho_{im}( \rho_{im}^{-1}(-)\bullet \rho_{im}^{-1}(-)):&
 H_*(\imm\ss) \times H_*(\imm\ss)\rightarrow \\ &H_{*-2n+1}(\imm  \ss)
\end{split}
\end{equation}

\bth
$H_{*+2n-1}(\imm (S^1,S^n))$ equipped with $\mu_{im}$ is a graded associative algebra.
\eth

After restricting to immersions with $k$ and $l$ double points, it reads

\begin{equation}
\begin{split}
\mu_{im}^{k,l}= \rho_{im}^{k+l}( {\rho^k_{im}}^{-1}(-)\bullet {\rho^l_{im}}^{-1}(-)):
 & H_*(\imm_{k}\ss) \times H_*(\imm_{l}\ss)\rightarrow \\ &H_{*-2n+1}(\imm_{k+l}  \ss)
\end{split}
\end{equation}

\begin{remark}  
We observe that all these theorems are corollaries of Gruher-Salvatore and Budney's work, we therefore don't propose a detailed proof. 
\end{remark}

\section{Comparison with Chas-Sullivan product}

So far we have only define the associative product on $H_{*}(\imm \ss)$ and not on $H_{*}(\im \ss)$. The algebra structure of the latter will be discussed in the last section  where we prove Theorem \ref{thm3} cited in the introduction. What is discussed in this section is the comparison of the algebra structure on $H_{*}(\imm \ss)$ with that of $H_{*}(LUS^n)$ equipped with Chas-Sullivan loop product.The  1-jet map  $\Psi: \gamma \mapsto (\gamma, \gamma'/\|\gamma'\|)$ is a map of fiberwise homotopy monoids making the following diagrams commute

\begin{equation}
\xymatrix{\imm \ss \ar[r]^-{\Psi}\ar[d] ^-{ev}& LUS^n \ar[d]^-{ev_{0}} \\  US^n \ar[r]^{id} & US^n},
\end{equation}
and 
\begin{equation}
\xymatrix{ SO(n+1)\times_{SO(n-1)} \im^l \rr \ar[r]^-{\Psi\circ \rho_{im}}\ar[d] & LUS^n \ar[d]^-{ev_{0}} \\  US^n \ar[r]^{id} & US^n}.
\end{equation}

Recall that  $ev(\gamma)= (\gamma(0), \gamma'(0)/\| \gamma'(0)\|)$ for $\gamma \in \emb \ss$ and $ev_{0}(\alpha)=\alpha(0)$ for $\alpha \in LUS^n $.

An immediate consequence of Proposition 11 in \cite{gs} is

\bth \label{thm-aux}
The map induced by 1-jet map $ (H_{*}(\imm \ss), \mu_{im}) \rightarrow (H_{*}(LUS^n), \bullet)$  is a map of graded associative algebras.  Here, $\bullet$ is the Chas-Sullivan loop product.
\eth

Using Hirsch-Smale theorem, one can state the theorem above differently. 

\begin{flushleft}
\textbf{Theorem} (Hirsch-Smale)
The 1-jet map $\im \ss \rightarrow L US^n$, given by $\gamma\mapsto  (\gamma, \gamma'/\| \gamma' \|)$, is a homotopy equivalence. 
\end{flushleft}

One can transfer the associative commutative product from $H_{*} (LUS^n)$ to $H_{*}(\im \ss)$, and then Theorem \ref{thm-aux} reads
 
\setcounter{section}{4} \setcounter{theorem}{1}
\bth
\label{thm3}
$H_{*+2n-1}(\im (S^1,S^n))$  is a graded associative algebra and the map induced by inclusion
$H_{*+2n-1}(\imm (S^1,S^n))\rightarrow H_{*+2n-1}(\im (S^1,S^n))$ is a map of algebras.
\eth
\section{Desingularization morphism}

We first introduce $\Q_{k}^n$ the \emph{desingularization spaces}
and the \emph{desingularization morphism} $\sigma_{k}: \Q_{k}^n\rightarrow \emb \rr$.  One should think of $\Q_{k}^n$  as the space of $k$-singular knots decorated with $k$ tangent vectors at the singular points, to be used to desingularize the knot (see \cite{CCL}). 

Let $G_{2,n}= \so/SO(2)\times SO(n-2)$ be the Grassmanian manifold of oriented 2-plane in $\R^n$. One has the canonical fibration $Q_{n} \rightarrow G_{2,n}$ whose fibers are $S^{n-3}$. Here $Q_{n}$ is the space of pairs $(\pi,v)$ where $\pi$ is a 2-plane and $v$ is  a vector perpendicular to $\pi$. More precisely,
$$
Q_{n}= \so \times_{(SO(2)\times SO(n-2))} S^{n-3}
$$
where the action of $SO(2)\times SO(n-2)$ is trivial for the first factor and for the second factor is given by the natural action of $SO(n-2)$ on $S^{n-3}\subset \R^{n-2}$. Then $Q_{n}^{\times k} $ fibers over $ G_{2,n}^{\times k}$ through the map $(\pi_{i}, v_{i}) _{1\leq i \leq k}\mapsto (\pi_{i})_{1\leq i \leq k}$.

Then we define $\Q_{n}^{k}$ to be the total space of the pull back of the fibration $Q_{n}^{\times k} \rightarrow G_{2,n}^{\times k}$  via the natural map $ r_{k}:\im^{l}_{k}\rr \rightarrow G_{2,n}^{\times k}$, given by sending an immersion with $k$ singular double points to the $k$ 2-planes define at the singularities by the tangent vectors. We have the following commutative diagram of fibrations,
\begin{equation} \label{comm-fibres}
\xymatrix{\Q_n^k \ar[r] \ar[d]&Q_{n}^{\times k} \ar[d] \\  \im^{l}_{k}\rr \ar[r]^-{r_{k}} & G_{2,n}^{\times k}}
\end{equation}
where the fibers are $({S^{n-3}})^{\times k}$.  

The \emph{desingularization morphism} $\sigma: \coprod _{k\geq 1} \Q_n^{k} \rightarrow \emb^l \rr$ is defined by resolving the singularities using the normal vectors given at those points. More explicitly,  let $(\gamma, (v_{i})_{1\leq i\leq k}) \in \Q_{k}$ where $\gamma: \R \rightarrow \R^n $ is a long immersion with double singular points $\gamma(t^1_{1})=\gamma(t^1_{2}), \gamma(t^2_{1})=\gamma(t^2_{2}),\cdots , \gamma(t_{1}^{k})=\gamma(t_{2}^{k})$ and $v_{1}, \cdots, v_{k}$ are the vectors normal to the tangent plane at the singularities. Consider
the bump functions,
\begin{equation}\label{reso}
\alpha_{i}(t)=
\begin{cases}
0 &\text{   if   } t\notin [t^{i}_{a_{i}}-\epsilon, t^{i}_{a_{i}}+\epsilon]\\
(-1)^{a_{i}}\delta \exp(\frac{-1}{(t-t^{i}_{a_{i}})^2 -\epsilon^2} ) v_i&  \text{ if }  t\in[t^{i}_{a_{i}}-\epsilon, t^{i}_{a_{i}}+\epsilon]
\end{cases}
\end{equation}
where each $a_{i}$ is 1 or 2.  Let the embedding  $\sigma((\gamma, (v_{i})_{ i=1}^{k})):\R\rightarrow \R^n$ be 
\begin{equation}
\sigma_k^n((\gamma, (v_{i})_{ i=1}^{k}))(t)=\gamma(t)+\sum_{i=1}^{k} \alpha_{i}(t).
\end{equation}
To define $\sigma_k^n$  one has to make a choice for $a_{i}$, but it turns out that these choices give rise to homologous cycles. Moreover, our construction obviously depends on good choices of $\epsilon$ and $\delta$. One has to make sure that each $t^i_{1}$ and $t^{i}_{2}$ are at least $\epsilon$ away from each other and $\delta$ is small enough not to create any new singularity. To resolve this problem, one considers the space ${\Q'}_{k}^n$ consisting of triples $ [ (\gamma, (v_{i})_{1\leq i\leq k}), \epsilon, \delta]$ for which the construction works. One has the fibration ${\Q'}^n_{k}\rightarrow \Q^n_{k}$, given by projecting on the first factor $(\gamma, (v_{i})_{1\leq i\leq k})$, whose fibers are contractible. Therefore ${\Q'}_{k}^n$ and $\Q_{k}^n$ have the same homotopy type.   Note that both ${\Q'}_{k}^n$ and $\Q_{k}^n$ thought of long immersions decorated with the some normal vectors at the singularities have a natural homotopy associative product.

Thus we have,

\begin{lemma}(See \cite{CCL}) Desingularization morphism $\sigma_{k}^n: \Q_{k}^n \rightarrow  \emb^l \rr $ is homologically well-defined and by the abuse of notation we have
$$
\sigma_{k}^n: H_{*}(\Q_{k}^n)\rightarrow H_{*}(\emb^l \rr).
$$
which is a map of associative algebras with respect to the product induced on the homologies by the concatenation products on $\Q_{k}^n$  and $\emb^l \rr$.
\end{lemma}

Moreover, $\sigma_k^n$ is clearly $SO(n-1)$ equivariant, therefore it induces a homologically well-defined map

$$
\Sigma_{k}^n: H_{*}(SO(n+1)\times_{SO(n-1)}\Q_{k}^n)\rightarrow H_{*}(SO(n+1)\times_{SO(n-1)} \emb^l \rr)  \\.
$$

We have a map

\begin{equation} 
\xymatrix{
\coprod _k  \Q_{k}^n \ar[r]^-{\coprod _k \sigma_{k}^n} \ar[d]&  \emb^l \rr
 \ar[d] 
\\
\coprod _k  SO(n+1)\times_{SO(n-1)}\Q_{k}^n \ar[r]^-{\coprod _k \Sigma_{k}^n} \ar[d]& SO(n+1)\times_{SO(n-1)} \emb^l \rr
 \ar[d] \\ 
SO(n+1)/SO(n-1)\ar[r]^-{id} &SO(n+1)/SO(n-1)}
\end{equation}
of fiberwise monoids, hence: 

\begin{lemma}  \label{lemm2}
$$
\oplus _k \Sigma_{k}^n: (\oplus _k H_{*}(SO(n+1)\times_{SO(n-1)}\Q_{k}^n),\ast)\rightarrow (H_{*}(SO(n+1)\times_{SO(n-1)} \emb^l \rr), \mu_{im})$$
is map of associative algebras.
\end{lemma}

We need another map in order to complete the construction of the desingularization morphism from the homology of the space of  immersions in $S^n$ to the homology of knot space.

\begin{proposition}(See \cite{BG})
Let $F\rightarrow E \rightarrow M$ be an oriented fibration whose fiber has the homotopy type of a closed manifold of dimension $p$.  Then there is a natural map
$f!: H_*(M)\rightarrow H_{*+p}(E)$. This is the dual of the map on cohomology given by integration over fibres.
\end{proposition}

A description using spectral sequence is the following composite

$$
E^2_{p,n}= H_p(M,\mathcal{H}_n(F))= H_p(M,\mathbb{Z}) \hookrightarrow  E^{\infty}_{p,n} \Rightarrow H_{p+n}(E)
$$

Here $\mathcal{H}_q(F)$ stands for the local coefficient system which is trivial in this case and has the fundamental class of $F$ as a generator (because the fibration is oriented). In this description one can see clearly the naturality property.

 Applying this construction to the oriented fibration $ \Q^{n}_{k} \rightarrow  \im_{k}^l \rr $  and
 $SO(n+1)\times_{SO(n-1)} \im_{k}^l \rr \rightarrow SO(n+1)\times_{SO(n-1)}  \Q_{k}^{n}$

\begin{proposition} The fibration $ \Q^{n}_{k} \rightarrow  \im_{k}^l \rr $ gives rises to maps
$$\theta_{n}^k: H_{*}( \im_{k}^l \rr )\rightarrow H_{*+k(n-3)}(\Q_{k}^{n})$$
and

$$
\Theta_{n}^k: H_{*}(SO(n+1)\times_{SO(n-1)} \im_{k}^l \rr )\rightarrow H_{*+k(n-3)}(SO(n+1)\times_{SO(n-1)}  \Q_{k}^{n})$$

\end{proposition}

It follows immediately from the spectral sequence  description of the Gysin maps
which can easily accomodates the multiplicative structure of the  homologies involved.

\begin{lemma} The  Gysin maps $\Theta_k^n$'s  are multiplicative in the sense that    

$$\Theta_k^n(x)\ast  \Theta_l^n(y)=\Theta_{k+l}^n(\mu_{im}(x,y)).$$

\end{lemma}

Then the composition
\begin{equation}
\begin{split}
\Sigma_{k}^n\circ \Theta^{n}_k:  &H_{*}(SO(n+1)\times_{SO(n-1)} \im_{k}^l \rr ) \rightarrow  \\ &H_{*-k(n-3)}(SO(n+1)\times_{SO(n-1)} \emb^l \rr)
\end{split}
\end{equation}
defines the desingularization morphism
 
\begin{equation}
\begin{split}
\phi^n_k=\Sigma_{k}^n\circ \Theta^{n}_k:  H_{*}(\im_{k} \ss ) \rightarrow H_{*-k(n-3)}( \emb\ss)
\end{split}
\end{equation}

\bth
The desingularization morphisms $\phi_k^n: H_*(\imm_k \ss) \rightarrow H_{*-k(n-3)}(\emb \ss)$, $k\geq 0$, are compatible with the products i.e.  for $x\in \H_*(\imm_k \ss)$ and $y\in H_*(\imm_{l}(\ss))$
$$
\mu_{em}(\phi_k^n (x), \phi_{l}^n(y))= \phi_{k+l}^n(\mu_{im}^{k,l}(x,y))
$$
\eth
\begin{proof}
For $x,y \in H_*(\imm_k \ss) $

\begin{equation}
\mu_{em}(\phi_k^n (x), \phi_{l}^n(y))= \mu_{em}(\Sigma_{k}^n\circ \Theta^{n}_k(x), \Sigma_{l}^n\circ \Theta^{n}_l(y))= \Sigma_{k+l}^n (\Theta^{n}_k(x)\ast \Theta^{n}_l(y))
\end{equation}
by Lemma \ref{lemm2}.  The latter is $\Sigma_{k+l}^n\circ \Theta_{k+l}^n(\mu_{im}^{k,l}(x,y))=\phi_{k+l}^n \mu_{im}^{k,l}(x,y))$.

\end{proof}

\section{Some computations}

Let us consider the canonical inclusion $i:\imm(S^1,S^n)\rightarrow
\im(S^1,S^n)$, the aim of this section is two-fold 
\\
first we determine the
morphism $i_*$ induced in singular homology  
\\
secondly we prove that the
desingularization morphism is non-trivial. 
\\
\\
{\bf Notations} Let $M$ be a d-dimensional closed oriented manifold, $\mathbb H_*(M)=H_{*+d}(M,\mathbb Z)$ is the regraded homology intersection algebra of $M$, $\mathbb H_*(UM)=H_{*+2d-1}(UM,\mathbb Z)$ is the regraded homology intersection algebra of $UM$ and
$\mathbb H_*(LM,\mathbb Z)=H_{*+d}(LM,\mathbb Z)$ is the Chas-Sullivan regraded loop homology intersection algebra of $LM$.

\subsection{Homology of immersion spaces}

\subsubsection{Immersions in $\mathbb R^{n+1}$.} As a baby example let us first consider the case of long immersions $\im^l(\mathbb R,\mathbb R^{n+1})$ and of singular knots $\im(S^1,\mathbb R^{n+1})$ into an Euclidean space. Thanks to the Hirsch-Smale's theorem one knows that the maps, $\gamma \mapsto D\gamma $ and $\gamma \mapsto (\gamma, D\gamma) $ give homotopy equivalences between these spaces and the spaces $\Omega  S^{n}$ and $L U\mathbb R^{n+1}$,  where $D(\gamma):=\frac{\dot{\gamma}}{\|\acute{\gamma}\|}$. 

\begin{equation}
\xymatrix{ \im^l(\mathbb R,\mathbb R^{n+1})\ar[r]^-{D}  \ar[d]&  \Omega  S^n \ar[d] \\ \im(S^1,\mathbb R^{n+1}) \ar[r]^-{(id,D)} \ar[rd]_-{D_{0}}&  LU \R^{n+1}=L\R ^{n+1}\times LS^n \ar[d]^-{ev _{0}}\\ & S^n}
\end{equation}

We identify $U\mathbb R^{n+1}$ together with $\mathbb R^{n+1}\times S^n$, and we consider the canonical projection 
$pr_2:U\mathbb R^{n+1}\rightarrow S^n$. Thus the spaces $\im^l(\mathbb R,\mathbb R^{n+1})$ and $\im(S^1,\mathbb R^{n+1})$  are homotopy equivalent to $\Omega S^n$ and $L S^n$.
\\
The concatenation of long immersions in 
$\im^l(\mathbb R,\mathbb R^{n+1})$ obviously corresponds to the composition of based loops in $\Omega S^n$, we recall that the singular homology of these spaces together with this multiplicative structure is isomorphic to the tensor algebra generated by an element $u$ of degree $n-1$. The evaluation map 
$$D_0:\im(S^1,\mathbb R^{n+1})\overset{D}{\rightarrow} LS^n\overset{ev_{0}}{\rightarrow}S^{n+1}$$
can be used to endow the singular homology of $\im(S^1,\mathbb R^{n+1})$ with a product. We consider the diagram

\begin{equation}
\xymatrix{ \im(S^1,\mathbb R^{n+1})\times \im(S^1,\mathbb R^{n+1}) &  \im(S^1,\mathbb R^{n+1})\times_{S^n} \im(S^1,\mathbb R^{n+1}) \ar[l]_-{j} \ar[d]^-{Tr}\\
& \im_{(0,\ldots,0)}(\infty,\mathbb R^{n+1})}
\end{equation}
where $\im(S^1,\mathbb R^{n+1})\times_{S^n} \im(S^1,\mathbb R^{n+1})$ is the space of immersions $(\gamma_1,\gamma_2)$ such that
$D_0(\gamma_1)=D_0(\gamma_2)$, the space $\im_{(0,\ldots,0)}(\infty,\mathbb R^{n+1})$ is the subspace of $\im(S^1,\mathbb R^{n+1})\times_{S^n} \im(S^1,\mathbb R^{n+1})$ such that $\gamma_1(0)=\gamma_2(0)=(0,\ldots,0)$. The map $Tr$ is a translation map, it is defined by using the unique translation $T_1$ that sends $\gamma_1(0)$ to $(0,\ldots, 0)$ and the unique translation $T_2$ that sends $\gamma_2(0)$ to $(0,\ldots,0)$, $Tr(\gamma_1,\gamma_2)=(T_1(\gamma_1),T_2(\gamma_2))$. We also have a composition map
$$comp:\im_{(0,\ldots,0)}(\infty,\mathbb R^{n+1})\rightarrow \im(S^1,\mathbb R^{n+1}).$$
\\
As $j$ is a finite codimensional embedding we can define a Gysin morphism $j_!$ and a product :
$$-\bullet-:=comp_*\circ Tr_*\circ j_!\circ -\times-.$$ 
This multiplicative structure is isomorphic to the Chas-Sullivan product on $ H_{*+n}(LS^n)$, the isomorphism is induced by the map $D$. 
\\
The algebra $H_{*+n}(LS^n)$ has been computed by Cohen, Jones and Yan \cite[Theorem 2]{CJY}. For a closed oriented manifold $M$ they consider the Serre spectral sequence associated to the evaluation fibration 
$$\Omega M\rightarrow LM \rightarrow M$$
enriched with a multiplicative structure coming from the loop product of $\mathbb H_*(LM)$ lifted at the chain level. To be more precise using the Pontryagin product on $H_{*}(\Omega M)$ and the intersection product on $\H_{*}(M)=H_{*+n}(M)$ the second page of this spectral sequence is (\cite[Theorem 1]{CJY})
$$E^2_{p,q}:=\mathbb H_p(M,H_q(\Omega M))\Rightarrow \mathbb H_{p+q}(LM):=H_{*+n}(LM).$$  
For the spheres, they get 
\[
\H_{*}(\im(S^{1},\mathbb{R}^{n+1}))\cong\mathbb{H}_{*}(LS^{n})\cong\left\{ \begin{array}{c}
\Lambda(a)\otimes\mathbb{Z}[u]\quad for\: n\: odd\\
(\Lambda(b)\otimes\mathbb{Z}[a,v])/(a^{2},ab,2av)\quad for\: n\: even\end{array}\right.\]
with $a\in \mathbb H_{-n}( LS^n)\cong E^{\infty}_{-n,0}$, $b\in \mathbb H_{-1}(LS^n)\cong E^{\infty}_{-n,n-1}$, 
$u\in\mathbb H_{n-1}(LS^n)\cong E^{\infty}_{0,n-1}$ and $v\in\mathbb H_{2n-2}(LS^n)\cong E^{\infty}_{0,2n-2}$. 

\subsubsection{Immersions in $S^n$.} The computation of $\mathbb H_*(\im(S^1,S^n))$ has been completed by Le Borgne and the second author in \cite{CL}. In fact Le Borgne in \cite{L} considers the loop fibration 
$$LF\rightarrow LE\rightarrow LB$$
associated to a smooth fiber bundle of oriented closed manifolds
$$F\rightarrow E\rightarrow B$$
then he enriches the Serre spectral sequence with a multplicative structure coming from the Chas-Sullivan loop product on $LB$, $LE$ and $LF$ :
$$E^2_{p,q}:=\mathbb H_p(LB,\mathbb H_q(LF))\rightarrow \mathbb H_{p+q}(LE).$$
In particular he applies this new technique to the fibration 
$$LS^{n-1}\rightarrow LUS^n\rightarrow LS^n$$
and proves that for $n$ odd, the spectral sequence collapses at the $E^2$-term. In the even case, he proves that the Serre spectral sequence of 
$$\Omega US^n\rightarrow LUS^n\rightarrow US^n$$
collapses at the $E^2$-term. In each case there are serious extension issues which are all solved in \cite{CL} (see the main theorem) using Morse theoretical techniques: the authors filter the immersion spaces by the energy functional they get a new multiplicative spectral sequence that converge to $\mathbb H_*(\im(S^1,S^n))$. We recall that we have the following isomorphisms of algebras : 
\\
{\bf Even case.} When $n$ is even,
$$\mathbb H_*(\im(S^1,S^n))\cong \mathbb H_*(US^n,H_*(\Omega US^n))$$
we also recall that 
$$\mathbb H_*(US^n)\cong  \Lambda(a_{-n+1},b_{-2n+1})/(2a_{-n+1},a_{-n+1}b_{-2n+1})$$
and that 
$$H_*(\Omega US^n)\cong \mathbb Z[u_{n-2},v_{2n-2}]/(2u_{n-2}).$$
The homology of $US^n$ can be computed using the Serre spectral sequence of the fibration 
$$S^{n-1}\rightarrow US^n\rightarrow S^n$$
and $H_*(\Omega US^n)$ using the Serre spectral sequence of
$$\Omega S^{n-1}\rightarrow \Omega US^n\rightarrow \Omega S^n.$$
{\bf Odd case.} For $n$ odd we have
$$\mathbb H_*(\im(S^1,S^n))\cong \mathbb H_*(LS^n)\otimes \mathbb H_*(LS^{n-1}).$$ 
\\
\\
In order to determine $i_*:H_{*}(\imm(S^1,S^n))\rightarrow H_{*}(\im(S^1,S^n))$, we need to consider the two fibrations 
$$\imm_{*}(S^1,S^n)\rightarrow \imm(S^1,S^n)\rightarrow S^n$$
$$\imm_{*,v}(S^1,S^n)\rightarrow \imm(S^1,S^n)\rightarrow US^n,$$
where $\imm_{*}$ denotes the space of based immersions and $\imm_{*,v}$ is its subspace consisting of those immersions with a given unit tangent vector at the base point. 

We notice that the second fibration admits a section $geod$ this map sends the unit tangent vector $v\in US^n$ to the unique great circle such that $(\gamma(0),\dot{\gamma}(0))=v$. 

\bth For $n$ even, 
$$\mathbb H_*(\imm(S^1,S^n))\cong\mathbb H_*(US^n,H_*(\Omega S^{n-1}))\cong \mathbb H_*(US^n)\otimes H_*(\Omega S^{n-1}).$$
The morphism $i_*$ is induced by the canonical morphism
$$\mathbb H_*(US^n,H_*(\Omega S^{n-1}))\rightarrow \mathbb H_*(US^n,H_*(\Omega US^{n}))\cong \mathbb H_*(\im(S^1,S^n))$$
which is the indentity on $\mathbb H_*(US^n)$ and the morphism 
$$j_*:H_*(\Omega S^{n-1})\rightarrow H_*(\Omega US^n)$$
given by the inclusion of a fibre $j:S^{n-1}\rightarrow US^n$.
For $n$ odd,
$$\mathbb H_*(\imm(S^1,S^n))\cong\mathbb H_*(S^n)\otimes \mathbb H_*(LS^{n-1}),$$
in this case $i_*$ is the injection
$$\mathbb H_*(S^n)\otimes \mathbb H_*(LS^{n-1})\stackrel{c_*\otimes Id}{\longrightarrow}\mathbb H_*(LS^n)\otimes \H_*(LS^{n-1}),$$
where $c:S^n\hookrightarrow LS^n$ is the inclusion of the constant loops.
\eth

\begin{proof}
({\bf the even case}). Let us consider the Serre spectral sequence associated to the fibration 
$$\imm_{*,v}(S^1,S^n)\rightarrow \imm(S^1,S^n)\rightarrow US^n.$$
The $E^2$ term of this spectral sequence is given by
$$E^2_{p,q}= \mathbb H_p(US^n,\mathbb H_q(\imm_{*,v}(S^1,S^n))).$$
Because $H_*(\Omega S^n)$ is torsion free we have the canonical isomorphism 
$$\mathbb H_p(US^n)\otimes H_q(\Omega S^{n-1}).$$
Let us recall that in the even case we have
$$\mathbb H_*(US^n)\cong \Lambda(a_{-n+1},b_{-2n+1})/(2a_{-n+1},a_{-n+1}b_{-2n+1})$$
and 
$$H_*(\Omega S^{n-1})\cong \mathbb Z[u_{n-2}].$$
In order to compute the other terms, we use the fact that this spectral sequence is multiplicative. As it has a section one has that the differentials
$d_k(a_{-n+1})=d_k(b_{-2n+1})=0$ for $k\geq 2$, and we also have $d_k(u_{n-2})=0$ for degree reasons. We conclude that this spectral sequence collapses at the $E^2$ term. Also for degree reasons there is no extensions issues thus we have that 
$$\mathbb H_*(\imm(S^1,S^n))\cong \Lambda(a_{-n+1},b_{-2n+1})/(2a_{-n+1},a_{-n+1}b_{-2n+1})\otimes \mathbb Z[u_{n-2}].$$  
\\
({\bf The odd case}). Let us consider the Serre-spectral sequence associated to the fibration 
$$\imm_{*}(S^1,S^n)\rightarrow \imm(S^1,S^n)\rightarrow S^n$$
as in the even case, for degree reasons this spectral sequence collapses at the $E_2$-term and there is no extensions issues.
\end{proof}  

\subsection{Non-triviality of the desingularization morphism} 

Let us prove that the
desingularization morphism 
$$\sigma_2:H_0(\im_2(S^1,S^n))\rightarrow H_{2n-6}(\emb(S^1,S^n))$$
is non-trivial. We suppose that $n>4$ and we follow R. Budney's construction,
in \cite{budney3} he gives a generator of $\pi_{2n-6}(\emb^l(\mathbb R,\mathbb
R^n))\cong \mathbb Z$ which is the first non-trivial homotopy group \cite[Proposition 3.9]{budney3}. 
\\
Let us describe this generator geometrically:  take a long immersion $\gamma$ together with two regular double points : we suppose that we
have $\gamma(t_1)=\gamma(t_3)$ and $\gamma(t_2)=\gamma(t_4)$ with
$t_1<t_2<t_3<t_4$, and that the resolution is the trefoil knot. 
Associated to this long immersion we have a desingularization map 
$$\text{des}:S^{n-3}\times S^{n-3}\rightarrow \emb^l(\mathbb R,\mathbb R^n).$$   
If $s_{2n-6}$ is the fundamental class of $S^{n-3}\times S^{n-3}$ then in homology
$\text{des}_*(s_{2n-6})$ is the generator of $H_{2n-6}(\emb^l(\mathbb R,\mathbb
R^n))$ (see \cite[Theorem 3.13]{budney3}).
\\
Let us take a closure $\gamma'\in \im_2(S^1,S^n)$ of $\gamma$ then\\
1) $H_{2n-6}(\emb(S^1,S^n))\cong \mathbb{Z}$ by using the Serre's spectral sequence of
the fibration $\emb(S^1,S^n)\rightarrow US^n$ and the fact that this spectral
has a section we have $H_{2n-6}(\emb(S^1,S^n))\cong E_2^{0,2n-6}$.   
\\
2) If $k_{2n-6}$ is a generator of $H_{2n-6}(\emb(S^1,S^n))$ then we have
$\sigma_2(1_{\text{tref}})=k_{2n-6}$ where $1_{\text{tref}}$ is the generator the $0$-th
homology group of the connected component containing $\gamma'$.

\end{document}